\documentclass{article}

\newtheorem{fed}{\textbf{Definition}}[section]
\newtheorem{thm}[fed]{\textbf{Theorem}}
\newtheorem{lemma}[fed]{\textbf{Lemma}}

\newtheorem{prop}[fed]{\textbf{Proposition}}

\usepackage{amssymb,bbm,graphicx,epsfig,psfrag,epic,eepic,latexsym}
\usepackage{amsmath}
\usepackage{mathrsfs}
\usepackage[dvips]{color}

\begin{document}
\title{A compactness theorem for frozen planets}
\author{Urs Frauenfelder}
\maketitle

\begin{abstract}
We study the moduli space of frozen planet orbits in the Helium atom for an interpolation between
instantaneous and mean interactions and show that this moduli space is compact. 
\end{abstract}

\section{Introduction}

Since the beginnings of quantum mechanics the dynamics of the Helium atom is intriguing people. 
Different from the Hydrogen atom the dynamics of the Helium atom is not completely integrable. In particular,
its phase space is not foliated by invariant tori and the EKB method cannot be applied to it. On the other
hand the Hamiltonian of Helium is invariant under simultanuous rotation of the two electrons or particle
interchange so that due to these symmetries period orbits are usually not isolated and Gutzwiller's trace
formula \cite{gutzwiller} cannot be directly applied neither. 
\\ \\
An interesting periodic orbit for the Helium problem was discovered numerically by Wintgen, Richter and Tanner see \cite{wintgen-richter-tanner, tanner-richter-rost}, which plays an important role in the semiclassical treatment of the Helium atom. For this orbit both electrons lay on a ray on the same side of the nucleus.
The inner electron collides with the nucleus and bounces back, while the outer electron (the frozen planet)
stays almost stationary simultanuously attracted by the nucleus and repelled by the inner electron. An interesting aspect of this periodic orbit is that to the authors knowledge it does not fit into a family
of periodic orbits starting from a completely integrable system. In fact if one ignores the interaction between
the two electrons one obtains a completely integrable system. But ignoring the interaction between the electrons both electrons are just attracted by the nucleus and therefore have to fall in it and there is nothing like a frozen planet orbit. 
\\ \\
In order to get a more tractable system the author replaced in \cite{frauenfelder} the instantaneous interaction of the two electrons by a mean interaction. It turns out that for mean interaction between the electrons the outer electron is really frozen, i.e.,\,it becomes stationary. One obtains a delay equation
for the inner electron and it is shown analytically in \cite{frauenfelder} that there exists a unique nondegenerate solution. 
\\ \\
In this paper we interpolate linearly between the instantaneous interaction and the mean interaction between
the two electrons. It is interesting to note that there is a special type of frozen planet orbit. At time
$t=0$ both the inner and the outer electron have vanishing velocity. The inner electron is than accelerated in
direction of the nucleus in which it falls after some time. Suppose now that at the instant where the
inner electron collides with the nucleus the outer electron has vanishing velocity again. Then one can let the movie run backwards. The inner electron jumps out of the nucleus and goes back to its initial position. Meanwhile the outer electron goes back as well to its initial condition so that at the end both electron
are back at their initial position again with zero velocity. One has a periodic orbit. We refer to such periodic orbits as \emph{symmetric frozen planets}. Here is the description of this moduli problem in formulas.
For a homotopy parameter $r\in [0,1]$ one looks at solutions
$$q_1 \in C^\infty\big([0,1],(0,\infty)\big),\quad q_2 \in C^0\big([0,1],[0,\infty)\big) \cap
C^\infty\big([0,1),(0,\infty)\big)$$
of the following moduli problem
\begin{equation}\label{modul}
\left\{\begin{array}{cc}
\ddot{q}_1(t)=-\frac{2}{q_1(t)^2}+\frac{r}{(\overline{q}_1-\overline{q}_2)^2}+\frac{1-r}{(q_1(t)-q_2(t))^2},
& t \in [0,1]\\
\ddot{q}_2(t)=-\frac{2}{q_2(t)^2}-\frac{r}{(\overline{q}_1-\overline{q}_2)^2}-\frac{1-r}{(q_1(t)-q_2(t))^2},
& t \in [0,1)\\
q_2(t)<q_1(t), & t \in [0,1]\\
\dot{q}_1(0)=\dot{q}_1(1)=\dot{q}_2(0)=q_2(1)=0.
\end{array}\right.
\end{equation}
The variable $q_1(t)$ describes the inner electron on the ray $(0,\infty)$ at time $t$ and the variable 
$q_2(t)$ the outer one. The nucleus lies at the origin with whom the inner electron collides at time
$t=1$. The positive numbers $\overline{q}_1$ and $\overline{q}_2$ denote the mean value of $q_1$ respectively
$q_2$, defined by
$$\overline{q}_1=\int_0^1 q_1(t)dt, \qquad \overline{q}_2=\int_0^1 q_2(t) dt.$$
For $r=0$ the first two equations become the second order ODE
$$\left\{\begin{array}{cc}
\ddot{q}_1(t)=-\frac{2}{q_1(t)^2}+\frac{1}{(q_1-q_2)^2},
& t \in [0,1]\\
\ddot{q}_2(t)=-\frac{2}{q_2(t)^2}-\frac{1}{(q_1-q_2)^2},
& t \in [0,1).
\end{array}\right.$$
The first term on the righthand side describes the attraction by the nucleus whose charge is two, since
the nucleus of the Helium atom consists of two protons. The second term describes the repulsion of the two
electrons. For $r=1$ the first two equations become
$$\left\{\begin{array}{cc}
\ddot{q}_1(t)=-\frac{2}{q_1(t)^2}+\frac{1}{(\overline{q}_1-\overline{q}_2)^2},
& t \in [0,1]\\
\ddot{q}_2(t)=-\frac{2}{q_2(t)^2}-\frac{1}{(\overline{q}_1-\overline{q}_2)^2},
& t \in [0,1).
\end{array}\right.$$
In this case the instantaneous interaction between the two electrons is replaced by the interaction of
their mean positions. This is not an ODE anymore, but an equation involving delay and prolay. 
\\ \\
The main result of this paper tells us that the moduli space of solutions of (\ref{modul}) is compact.
There are two horror scenarios which have to be ruled out. The first horror scenario is ionization meaning
that the outer electron escapes to infinity. The other horror scenario is that the two electrons come arbitrary
close together. There are two ways how the later scenario could occur. Namely both electrons fall
simultaneously in the nucleus which leads to a triple collision, or the energy which is not bounded a priori
explodes. The following is the main result of this paper which rules out these horror scenarios.  
 
\begin{thm}\label{main}
There exists a constant $\kappa$ such that for every solution $(q_1,q_2,r)$ of the problem (\ref{modul}) one
has
$$\max_{t \in [0,1]}\bigg\{q_1(t), \frac{1}{q_1(t)-q_2(t)}\bigg\} \leq \kappa.$$ 
\end{thm}
In view of the theorem, given a sequence of solutions $(q_1^\nu,q_2^\nu,r^\nu)$ of problem (\ref{modul})
one can find a convergent subsequence which then by a usual bootstrapping argument is again a solution
of problem (\ref{modul}). The author does not know if nonsymmetric frozen planet orbits actually exist. If they exist it is an open question if the compactness result can be extended to the nonsymmetric case.
\\ \\
In a joint work with K.\,Cieliebak and E.\,Volkov the author is currently studying a variational approach
to frozen planet orbits \cite{cieliebak-frauenfelder-volkov}. Theorem~\ref{main} together with the 
variational approach then leads to a well-defined Euler characteristic at least for symmetric frozen planet
orbits. A more difficult question is, if a homology can be defined whose chain complex is generated by
symmetric frozen planet orbits. To define such a homology one needs to generalize the compactness result to gradient flow lines.   
\\ \\
In joint work with P.\,Albers, F.\,Schlenk, and J.\,Weber the author started to generalize Floer homology
to Hamiltonian delay equations \cite{albers-frauenfelder-schlenk1, albers-frauenfelder-schlenk2, albers-frauenfelder-schlenk3, frauenfelder-weber}. The author believes that the frozen planet problem can trigger a lot of research in this direction and can become an important testing ground how far Floer homology can be further developed.
\\ \\
\emph{Acknowledgements: } The author acknowledges partial support by DFG grant FR 2637/2-2.

\section{Proof of the main result}

We assume that $(q_1,q_2)$ is a solution of problem (\ref{modul}) for some $r \in [0,1]$.
Note that $\ddot{q}_2$ is always negative. Therefore in view of its initial condition
we have $\dot{q}_2(t)<0$ for every $t \in (0,1)$, so that $q_2$ is strictly monoton decreasing. 
The following lemma tells us that in contrast to $q_2$ the variable $q_1$ is monotone increasing. 
\begin{lemma}\label{pos}
If $r<1$, then for every $t \in (0,1)$ we have $\dot{q}_1(t)>0$. In particular, $q_1$ is strictly monotone increasing. If $r=1$, then $q_1$ is constant.
\end{lemma}
\textbf{Proof: } That $q_1$ is constant in the case $r=1$ is proved in \cite[Lemma\,3.1]{frauenfelder}. It
suffices therefore to consider the case $r<1$. We prove it in two steps.
\\ \\
\textbf{Step\,1: } For every $t_0 \in [0,1)$ such that $\dot{q}_1(t_0)=0$ we have $\ddot{q}_1(t_0)>0$.
\\ \\
We first note   
\begin{eqnarray*}
\ddot{q}_1-\ddot{q}_2&=&\frac{2}{q_2^2}-\frac{2}{q_1^2}+\frac{2r}{(\overline{q}_1-\overline{q}_2)^2}+\frac{2(1-r)}{(q_1-q_2)^2}\\
&=&\frac{2(q_1^2-q_2^2)}{q_1^2 q_2^2}+\frac{2r}{(\overline{q}_1-\overline{q}_2)^2}+
\frac{2(1-r)}{(q_1-q_2)^2}\\
&>&0.
\end{eqnarray*}
In view of the initial condition this implies
\begin{equation}\label{mon}
\dot{q}_1(t)-\dot{q}_2(t) > 0, \quad t>0.
\end{equation}
The jerk of $q_1$ is given by 
\begin{equation}\label{jerk}
\dddot{q}_1=\frac{4\dot{q}_1}{q_1^3}-\frac{2(1-r)(\dot{q}_1-\dot{q}_2)}{(q_1-q_2)^3}
\end{equation}
In view of the initial conditions we have
$$\dddot{q}_1(0)=0$$
and in view of (\ref{mon}) and $r<1$ it holds that
\begin{equation}\label{jerk2}
\dddot{q}_1(t)<\frac{4\dot{q}_1(t)}{q_1^3}, \quad t>0.
\end{equation}
We assume by contradiction that there exists
$$t_0 \in (0,1)$$
satisfying
\begin{equation}\label{t0}
\dot{q}_1(t_0)=0, \quad \ddot{q}_1(t_0) \leq 0.
\end{equation}
We define
$$t_1:=\inf\big\{t \in (t_0,1]: \dot{q}_1(t) \geq 0\big\}.$$
Here we use the convention, that if there is no $t \in (t_0,1]$ satisfying $\dot{q}_1(t) \geq 0$,
then $t_1=1$. 
In view of (\ref{jerk}) it follows that
$$\dddot{q}_1(t_0)=\frac{\dot{q}_2(0)}{(q_1(t_0)-q_2(t_0))^3}<0.$$
Therefore 
$$t_1>t_0.$$
In view of the definition of $t_1$ and (\ref{jerk2}) we have
$$\dddot{q}_1(t)<0, \quad t \in (t_0,t_1).$$
In view of (\ref{t0}) this implies that
$$\ddot{q}_1(t)<0, \quad t \in (t_0,t_1).$$
Combining this once more with (\ref{t0}) we conclude that
$$\dot{q}_1(t_1)<0.$$
By definition of $t_1$ this implies that
$$t_1=1$$
and therefore
$$\dot{q}_1(1)<0.$$
This contradicts the boundary condition in (\ref{modul}) and Step\,1 is proved. 
\\ \\
\textbf{Step\,2: } We prove the lemma.
\\ \\
By Step\,1 we conclude that $q_1$ cannot have a local maximum in $[0,1)$ but all its critical points
are strict local minima. By the boundary condition in (\ref{modul}) the function $q_1$ has a critical point
at time $0$ which therefore has to be a strict local minimum. Since there are no local maxima, there cannot
be any addition critical points of $q_1$ in $(0,1)$. In particular, the derivative $\dot{q}_1$ cannot change
sign. Since $q_1$ has a strict local minimum at time $0$ the sign of $\dot{q}_1(t)$ is positive for
$t$ close to $0$ and consequently is positive always. This finishes the proof of the lemma. \hfill $\square$
\begin{lemma}
For the starting point of $q_2$ we have the following lower bound
\begin{equation}\label{start}
q_2(0) \geq 1.
\end{equation}
\end{lemma}
\textbf{Proof: } By (\ref{modul}) we have for every $t \in [0,1)$ the estimate
$$\ddot{q}_2(t) \leq -\frac{2}{q_2(t)^2} \leq -\frac{2}{q_2(0)^2}$$
where the second inequality holds since $q_2$ is monotone decreasing. Using that  $\dot{q}_2(0)=0$ we obtain
from this the estimate
$$q_2(t) \leq -\frac{t^2}{q_2(0)^2}+q_2(0).$$
Taking the limit $t$ going to $1$ we get
$$0 \leq -\frac{1}{q_2(0)^2}+q_2(0)$$
implying
$$1 \leq q_2(0)^3$$
and therefore
$$1 \leq q_2(0).$$
This finishes the proof of the lemma. \hfill $\square$
\\ \\
We abbreviate
$$\Delta:=q_1(0)-q_2(0)$$
the distance between $q_1$ and $q_2$ at time $t=0$. Since the variable $q_2$ is decreasing and by Lemma~\ref{mon} the variable $q_1$ is increasing we have for every $t \in [0,1]$
\begin{equation}\label{distance}
q_1(t)-q_2(t) \geq \Delta.
\end{equation}
The following lemma tells us that $q_1$ and $q_2$ cannot come to close to each other.
\begin{lemma}\label{delr}
There exists a constant $c_1>0$ such that 
\begin{equation}\label{delest}
\Delta \geq (1-r)c_1. 
\end{equation}
\end{lemma}
\textbf{Proof: } We suppose that
$$\Delta \leq \frac{1}{2}$$
which we can assume without loss of generality after maybe shrinking the constant $c$. 
We define $t_0 \in (0,1)$ by the requirement
$$q_2(t_0)=q_2(0)-\Delta.$$
Note that
$$\ddot{q}_1+\ddot{q}_2=-\frac{2}{q_1^2}-\frac{2}{q_2^2}\geq -\frac{4}{q_2^2}$$
so that for $t \in [0,t_0]$ we have the estimate
$$\ddot{q}_1(t)+\ddot{q}_2(t) \geq -\frac{4}{(q_2(0)-\Delta)^2}\geq-\frac{4}{(1-\Delta)^2} \geq -16.$$
Since the velocity of both variables $q_1$ and $q_2$ at time $t=0$ vanish we get from that the estimate
$$\dot{q}_1(t_0)+\dot{q}_2(t_0) \geq -16t_0\geq-16$$
which we can rearrange to
$$-\dot{q}_2(t_0) \leq 16+\dot{q}_1(t_0).$$
Since $\dot{q}_2$ is negativ but by Lemma~\ref{pos} we have that $\dot{q}_1$ is positive we obtain from that the inequality
\begin{equation}\label{kin}
\dot{q}_2(t_0)^2 \leq 512+2\dot{q}_1(t_0)^2.
\end{equation}
Since $\ddot{q}_1+\ddot{q}_2<0$ we further have the estimate
$$q_1(t_0)-q_1(0) \leq q_2(0)-q_2(t_0)=\Delta$$
so that we obtain
\begin{equation}\label{dist}
q_1(t_0)-q_2(t_0)=\big(q_1(t_0)-q_1(0)\big)+\big(q_1(0)-q_2(0)\big)+\big(q_2(0)-q_2(t_0)\big) \leq 3\Delta
\leq \frac{3}{2}.
\end{equation}
We have the following preserved quantity
$$E=\frac{1}{2}\Big(\dot{q}_1(t)^2+\dot{q}_2(t)^2\Big)-\frac{2}{q_1(t)}-\frac{2}{q_2(t)}
+\frac{r\big(q_1(t)-q_2(t)\big)}{(\overline{q}_1-\overline{q}_2)^2}+\frac{1-r}{q_1(t)-q_2(t)}.$$
At time $t=0$ this computes to be
$$E=-\frac{2}{q_2(0)+\Delta}-\frac{2}{q_2(0)}+\frac{r\Delta}{(\overline{q}_1-\overline{q}_2)^2}
+\frac{1-r}{\Delta}$$
which we estimate using (\ref{start})
\begin{equation}\label{E1}
E \geq -4+\frac{1-r}{\Delta}.
\end{equation}
Note that there exists a constant $0<\varepsilon<1$ such that
$$\overline{q}_2 \leq (1-\varepsilon)q_2(0).$$
By the computations in Appendix~\ref{mean} the constant $\varepsilon$ can be chosen around $\tfrac{1}{4}$. 
Hence using Lemma~\ref{mon} and (\ref{start}) we estimate
$$\frac{1}{\overline{q}_1-\overline{q}_2} \geq \frac{1}{q_1(0)-\overline{q}_2} \geq
\frac{1}{q_2(0)-\overline{q}_2}\geq \frac{1}{\varepsilon q_2(0)} \geq \frac{1}{\varepsilon}.$$
Evaluating the preserved quantity $E$ at time $t=t_0$ we obtain from that as well as (\ref{kin}) and (\ref{dist})
\begin{eqnarray}\label{E2}
E&=&\frac{1}{2}\Big(\dot{q}_1(t_0)^2+\dot{q}_2(t_0)^2\Big)-\frac{2}{q_1(t_0)}-\frac{2}{q_2(t_0)}
+\frac{r\big(q_1(t_0)-q_2(t_0)\big)}{(\overline{q}_1-\overline{q}_2)^2}\\ \nonumber
& &+\frac{1-r}{q_1(t_0)-q_2(t_0)}\\ \nonumber
&\leq& \frac{3}{2}\dot{q}_1^2(t_0)+256+\frac{3}{2\varepsilon^2}+\frac{1-r}{q_1(0)-q_2(0)+\Delta}\\ \nonumber
&=& \frac{3}{2}\dot{q}_1^2(t_0)+256+\frac{3}{2\varepsilon^2}+\frac{1-r}{2\Delta}
\end{eqnarray}
Combining (\ref{E1}) and (\ref{E2}) we get the inequality
$$\frac{1-r}{2\Delta}\leq \frac{3}{2}\dot{q}_1^2(t_0)+260+\frac{3}{2\varepsilon^2}.$$
We abbreviate
$$c_0:=260+\frac{3}{2\varepsilon^2}.$$
Since $\dot{q}_1(t_0)$ is positive by Lemma~\ref{mon} we obtain the estimate
\begin{equation}\label{velest}
\dot{q}_1(t_0) \geq \frac{\sqrt{1-r}}{\sqrt{3\Delta}}-c_0.
\end{equation}
Using Lemma~\ref{mon} and (\ref{start}) again the acceleration of $q_1$ is estimated from below for every $t \in [0,1]$ by
\begin{equation}\label{ac1}
\ddot{q}_1(t) \geq -\frac{2}{q_1(t)^2} \geq -\frac{2}{q_1(0)^2} \geq -\frac{2}{q_2(0)^2}\geq -2.
\end{equation}
Since the velocity of $q_1$ at time $t=1$ vanishes we obtain from this combined with (\ref{velest})
\begin{eqnarray*}
0&=&\dot{q}_1(1)\\
&=&\dot{q}_1(t_0)+\int_{t_0}^1 \ddot{q}_1(t)dt\\
&\geq& \frac{\sqrt{1-r}}{\sqrt{3\Delta}}-c_0-2(1-t_0)\\
&\geq& \frac{\sqrt{1-r}}{\sqrt{3\Delta}}-c_0-2
\end{eqnarray*}
implying
$$\Delta \geq \frac{1}{3(c_0+2)^2}(1-r)$$
so that (\ref{delest}) follows with
$$c_1:=\frac{1}{3(c_0+2)^2}.$$
This finishes the proof of the lemma. \hfill $\square$
\\ \\
We would like to replace the estimate for $\Delta$ in Lemma~\ref{delr} by a uniform one not depending
on $r$. For that purpose we need the following result on the average positions of $q_1$ and $q_2$.
\begin{lemma}\label{meanest}
The mean values of $q_1$ and $q_2$ satisfy the following inequality
\begin{equation}\label{meane}
\overline{q}_2 \leq \bigg(1-\sqrt{\frac{r}{1+r}}\bigg)\overline{q}_1.
\end{equation}
\end{lemma}
\textbf{Proof: }From Lemma~\ref{mon} it follows that
$$\overline{q}_1 \leq q_1(1)$$
and $q_1$ attains at time $t=1$ a local maximum so that
$$\ddot{q}_1(1) \leq 0.$$
Combining these inequalities with (\ref{modul}) we estimate
\begin{eqnarray*}
0 &\geq&\ddot{q}_1(1)\\
&=&-\frac{2}{q_1(t)^2}+\frac{r}{(\overline{q}_1-\overline{q}_2)^2}+\frac{1-r}{(q_1(t)-q_2(t))^2}\\
&=&-\frac{2}{q_1(t)^2}+\frac{r}{(\overline{q}_1-\overline{q}_2)^2}+\frac{1-r}{q_1(t)^2}\\
&=&-\frac{1+r}{q_1(t)^2}+\frac{r}{(\overline{q}_1-\overline{q}_2)^2}\\
&\geq&-\frac{1+r}{\overline{q}_1^2}+\frac{r}{(\overline{q}_1-\overline{q}_2)^2}
\end{eqnarray*}
which we rewrite
$$\frac{1+r}{\overline{q}_1^2}\geq \frac{r}{(\overline{q}_1-\overline{q}_2)^2}$$
or equivalently
$$\frac{\overline{q}_1^2}{1+r}\leq \frac{(\overline{q}_1-\overline{q}_2)^2}{r}.$$
Taking square roots we obtain
$$\overline{q}_1 \leq \sqrt{\frac{1+r}{r}}(\overline{q}_1-\overline{q}_2)$$
so that
$$\sqrt{\frac{1+r}{r}}\overline{q}_2 \leq \bigg(\sqrt{\frac{1+r}{r}}-1\bigg)\overline{q}_1$$
from which (\ref{meane}) follows. \hfill $\square$
\\ \\
Now we are in position to improve Lemma~\ref{delor} with a uniform estimate from below for $\Delta$
not depending on $r$ anymore. 
\begin{prop}\label{delor}
There exists a constant $c_2>0$ such that $\Delta \geq c_2$.
\end{prop}
\textbf{Proof: } Using (\ref{ac1}) and $\dot{q}_1(1)=0$ we obtain the estimate
$$\dot{q}_1(1-t) \leq 2t$$
and therefore using integration by parts
\begin{eqnarray*}
\overline{q}_1&=&\int_0^1 q_1(t)dt\\
&=&\int_0^1 q_1(1-t)dt\\
&=&q_1(0)+\int_0^1 \dot{q}_1(1-t) t dt\\
&\leq&q_1(0)+\int_0^1 2t^2dt\\
&=&q_1(0)+\frac{2}{3}
\end{eqnarray*}
so that
\begin{equation}\label{meest1}
q_1(0) \geq \overline{q}_1-\frac{2}{3}.
\end{equation}
Since $\ddot{q}_2(t) \leq 0$ for every $t \in [0,1)$ and $\dot{q}_2(0)=0$ we have
\begin{equation}\label{meest2}
q_2(0) \leq 2 \overline{q}_2
\end{equation}
so that combined with (\ref{start}) we get the estimate
\begin{equation}\label{meest3}
\overline{q}_2 \geq \frac{1}{2}.
\end{equation}
From (\ref{meane}) we obtain for every $r \in [0,1]$ the inequality
\begin{equation}\label{meest4}
\overline{q}_1 \geq \frac{1}{1-\sqrt{\frac{r}{1+r}}}\,\overline{q}_2=\frac{\sqrt{1+r}}{\sqrt{1+r}-\sqrt{r}}
\,\overline{q}_2.
\end{equation}
Using (\ref{meest1}), (\ref{meest2}), and (\ref{meest4}) we estimate
\begin{eqnarray}\label{meest5}
\Delta&=& q_1(0)-q_2(0)\\ \nonumber
&\geq&\overline{q}_1-\frac{2}{3}-2\overline{q}_2\\ \nonumber
&\geq&\bigg(\frac{\sqrt{1+r}}{\sqrt{1+r}-\sqrt{r}}-2\bigg)\overline{q}_2-\frac{2}{3}.
\end{eqnarray}
Note that the function
$$f \colon [0,1] \to (0,\infty), \quad r \mapsto \frac{\sqrt{1+r}}{\sqrt{1+r}-\sqrt{r}}
=\frac{1}{1-\sqrt{\frac{r}{1+r}}}$$
is monotone increasing, since $r \mapsto \tfrac{r}{1+r}$ is monotone increasing, and
$$f\big(\tfrac{1}{3}\big)=\frac{1}{1-\sqrt{\frac{1}{4}}}=\frac{1}{1-\frac{1}{2}}=2$$
so that it follows from (\ref{meest3}) and (\ref{meest5}) that
\begin{equation}\label{meest6}
\Delta \geq \bigg(\frac{\sqrt{1+r}}{\sqrt{1+r}-\sqrt{r}}-2\bigg)\frac{1}{2}-\frac{2}{3}
=\frac{\sqrt{1+r}}{2\big(\sqrt{1+r}-\sqrt{r}\big)}-\frac{5}{3},\qquad r \geq \frac{1}{3}.
\end{equation}
Note that
$$f\big(\tfrac{144}{145}\big)=\frac{1}{1-\sqrt{\frac{144}{289}}}=\frac{1}{1-\frac{12}{17}}=\tfrac{17}{5}$$
so that since $f$ is monotone increasing
$$f(r)\geq \frac{17}{5}, \qquad r \geq \frac{144}{145}.$$
Therefore it follows from (\ref{meest6}) that
$$\Delta \geq \frac{17}{10}-\frac{5}{3}=\frac{51-50}{30}=\frac{1}{30}, \qquad r \geq \frac{144}{145}.$$
Combining this estimate with Lemma~\ref{delr} proves the proposition. \hfill $\square$
\\ \\
Our next goal is to derive an upper bound for the outer electron. We first start with an upper bound
for the inner electron.
\begin{lemma}\label{startup}
There exists a constant $c_3$ such that $q_2(0) \leq c_3$.
\end{lemma}
\textbf{Proof: } We recall (\ref{distance}) which tells us that
$$q_1(t)-q_2(t) \geq \Delta, \quad t \in [0,1],$$
which is an immediate consequence of Lemma~\ref{pos}.
In particular we have 
$$\overline{q}_1-\overline{q}_2 \geq \Delta.$$
We infer from (\ref{modul}) that
$$\ddot{q}_2(t) \geq -\frac{2}{q_2(t)^2}-\frac{1}{\Delta^2}, \quad t \in [0,1].$$
With Proposition~\ref{delor} it follows that
$$\ddot{q}_2(t) \geq -\frac{2}{q_2(t)^2}-\frac{1}{c_2^2}, \quad t \in [0,1].$$
Let $t_0 \in [0,1)$ be the time such that
$$q_2(t_0)=1.$$
Because $q_2$ is motonote decreasing we have
$$\ddot{q}_2(t) \geq -2-\frac{1}{c_2^2}, \quad t \in [t_0,1].$$
Since $\dot{q}_2(0)=0$ we deduce from that
$$q_2(0) \leq q_2(t_0)+\frac{t_0^2}{2}\bigg(2+\frac{1}{c_2^2}\bigg)
\leq 1+\frac{1}{2}\bigg(2+\frac{1}{c_2^2}\bigg)=2+\frac{1}{2c_2^2}.$$
Setting
$$c_3:=2+\frac{1}{2c_2^2}$$
the lemma follows. \hfill $\square$
\\ \\
Now we are in position to obtain an upper bound for the outer electron.
\begin{prop}\label{ion}
There exists a constant $c_4$ such that $q_1(1) \leq c_4$.
\end{prop}
\textbf{Proof: } Since $q_2$ is monotone decreasing we obtain from Lemma~\ref{startup} that
$$q_2(t) \leq c_3, \quad t \in [0,1]$$
and in particular
$$\overline{q}_2 \leq c_3.$$
Recall the inequality (\ref{ac1}) telling us that
$$\ddot{q}_2(t) \geq -2, \quad t \in [0,1].$$
Combining this with the fact that $q_1$ is monotone increasing by Lemma~\ref{pos}
we obtain the estimate
$$q_1(1) \geq q_1(t) \geq q_1(1)-1, \quad t \in [0,1]$$
and especially
$$\overline{q}_1 \geq q_1(1)-1.$$
Together with (\ref{modul}) these estimates imply
\begin{eqnarray*}
\ddot{q}_2(t) &\leq& -\frac{2}{q_1(1)^2}+\frac{r}{(q_1(1)-1-c_3)^2}+\frac{1-r}{(q_1(1)-1-c_3)^2}\\
&=&-\frac{2}{q_1(1)^2}+\frac{1}{(q_1(1)-1-c_3)^2}, \quad t \in [0,1].
\end{eqnarray*}
Since $\dot{q}_1(0)=\dot{q}_1(1)=0$ there exists $t_0$ satisfying $\ddot{q}_1(t_0)=0$ so that we
obtain the inequality
$$0 \leq -\frac{2}{q_1(1)^2}+\frac{1}{(q_1(1)-1-c_3)^2}$$
implying
$$q_1(1)^2 \geq 2\big(q_1(1)-1-c_3\big)^2.$$
Taking square roots on both sides we obtain the inequality
$$q_1(1) \geq \sqrt{2}\big(q_1(1)-1-c_3\big)$$
implying that
$$c_3+1 \geq \big(\sqrt{2}-1\big)q_1(1)$$
so that
$$q_1(1) \leq \frac{c_3+1}{\sqrt{2}-1}.$$
Hence setting
$$c_4:=\frac{c_3+1}{\sqrt{2}-1}$$
the proposition follows. \hfill $\square$ 
\\ \\
\textbf{Proof of Theorem~\ref{main}: } By Lemma~\ref{pos} we know that $q_1(t)$ is monotone increasing so
that combined with Proposition~\ref{ion} we have the estimate
\begin{equation}\label{ma1}
q_1(t) \leq c_4, \quad t \in [0,1].
\end{equation}
From (\ref{distance}) and Proposition~\ref{delor} we infer that
$$q_1(t)-q_2(t) \geq c_2, \quad t \in [0,1]$$
so that 
\begin{equation}\label{ma2}
\frac{1}{q_1(t)-q_2(t)} \leq \frac{1}{c_2}.
\end{equation}
Setting
$$\kappa:=\max\bigg\{c_4, \frac{1}{c_2}\bigg\}$$
the theorem follows from inequalities (\ref{ma1}) and (\ref{ma2}). \hfill $\square$

\appendix

\section{The average mean fall}\label{mean}

The interior electron of the frozen planet problem is in the free fall. In this appendix we discuss
mean values of some free falls. This is not really needed for the proof of the main result. On the other
hand a careful analysis of the mean free fall could be used to determine some constants occuring in
the proof more precisely and might be of use when trying to establish a homology theory which also involves
compactness results for gradient flow lines. Moreover, the mean free fall has its own mathematical beauty.
\\ \\
If the acceleration is constant $g>0$ then the free fall starting at height $q_0>0$ is the solution of the
initial value problem
$$q(0)=q_0,\quad \dot{q}_0=0, \quad \ddot{q}(t)=-g$$
whose explicit solution is given by
$$q(t)=q_0-\frac{g}{2}t^2.$$
If $\tau$ is the time of the free fall implicitly defined by
$$q(\tau)=0$$
then from the above formula one obtains explicitly
$$\tau=\sqrt{\frac{2q_0}{g}}.$$
The average position is then given by
\begin{eqnarray*}
\overline{q}&=&\frac{1}{\tau}\int_0^\tau q(t) dt=\frac{1}{\tau}\int_0^\tau\bigg(q_0-\frac{g}{2}t^2\bigg)dt
=\frac{1}{\tau}\bigg(q_0 \tau -\frac{g}{6}\tau^3\bigg)=q_0-\frac{g}{6}\cdot \frac{2q_0}{g}\\
&=&\frac{2}{3}q_0,
\end{eqnarray*}
so that the ratio of the average position and the initial position is given by
$$\kappa:=\frac{\overline{q}}{q_0}=\frac{2}{3}.$$
More generally, if the acceleration is given by the derivative $f'$ of a potential $f$, the free fall
is a solution of the initial value problem
$$q(0)=q_0, \quad \dot{q}_0=0, \quad \ddot{q}(t)=f'(q(t)).$$
One has the preserved quantity
$$\frac{1}{2}\dot{q}^2(t)+f(q(t))=f(q_0).$$
We assume that $f'<0$ so that $q$ is strictly decreasing.
Hence the velocity at time $t$ is given by
$$\dot{q}(t)=-\sqrt{2\big(f(q_0)-f(q(t)\big)}.$$
Using this formula the time of the free fall
can be computed as
$$\tau=\int_0^\tau dt=\int_0^{q_0}\frac{1}{\sqrt{2(f(q_0)-f(q))}}dq,$$
and the average position
\begin{eqnarray*}
\overline{q}=\frac{1}{\tau}\int_0^\tau q(t)dt=\frac{\int_0^{q_0}\frac{q}{\sqrt{2(f(q_0)-f(q))}}dq}{\int_0^{q_0}\frac{1}{\sqrt{2(f(q_0)-f(q))}}dq}=\frac{\int_0^{q_0}\frac{q}{\sqrt{f(q_0)-f(q)}}dq}{\int_0^{q_0}\frac{1}{\sqrt{f(q_0)-f(q)}}dq}.
\end{eqnarray*} 
Note that this expression is invariant under scaling the potential $f$ to $\mu f$ for $\mu>0$. We now
want to compute this for the homogeneous potentials 
$$f_\alpha(q)=-\frac{1}{q^\alpha}$$
for $\alpha>0$. 
We abbreviate by
$$\kappa(\alpha):=\frac{\overline{q}}{q_0}$$
the ratio between average position and initial position for the free fall with respect to the 
potential $\alpha$. We have the following proposition
\begin{prop}\label{free}
 The ratio $\kappa(\alpha)$ is given by
$$\kappa(\alpha):=\frac{\Gamma\big(\frac{4+\alpha}{2\alpha}\big)\Gamma\big(\frac{1+\alpha}{\alpha}\big)}{
\Gamma\big(\frac{2+\alpha}{\alpha}\big)\Gamma\big(\frac{2+\alpha}{2\alpha}\big)}.$$
\end{prop}
Let us look at some special values. Using $\Gamma(x+1)=x\Gamma(x)$ and $\Gamma(n+1)=n!$ we obtain for
the Newtonian potential $-\tfrac{1}{q}$, i.e., $\alpha=1$,
$$\kappa(1)=\frac{\Gamma\big(\frac{5}{2}\big)\Gamma(2)}{\Gamma(3)\Gamma\big(\frac{3}{2}\big)}
=\frac{\frac{3}{2}\Gamma\big(\frac{3}{2}\big)}{2\Gamma\big(\frac{3}{2}\big)}=\frac{3}{4}.$$
This is a bit bigger than $\tfrac{2}{3}$ one obtains for the free fall with constant acceleration. 
For the potential $-\tfrac{1}{q^2}$, i.e., $\alpha=2$, one obtains using $\Gamma\big(\tfrac{1}{2}\big)=\sqrt{\pi}$ the transzendental number
$$\kappa(2)=\frac{\Gamma\big(\frac{3}{2}\big)\Gamma\big(\frac{3}{2}\big)}{\Gamma(2)\Gamma(1)}=
\frac{\Gamma\big(\frac{1}{2}\big)^2}{4}=\frac{\pi}{4}.$$ 
\textbf{Proof of Proposition~\ref{free}: } Since the average position does not change if we
scale the potential by a positive factor we work with the potential $-\tfrac{1}{2q^\alpha}$ to avoid
annoying factors $\sqrt{2}$ in the computation. We have than
$$\dot{q}=-\sqrt{\frac{1}{q^\alpha}-\frac{1}{q_0^\alpha}}=-\frac{\sqrt{q_0^\alpha-q^\alpha}}
{(q_0q)^{\frac{\alpha}{2}}}$$
so that we obtain for the time of the free fall
\begin{eqnarray*}
\tau&=&q_0^{\frac{\alpha}{2}}\int_0^{q_0}\frac{q^{\frac{\alpha}{2}}}{\sqrt{q_0^\alpha-q^\alpha}}dq
\end{eqnarray*}
Changing variables
$$q=q_0 (\cos \theta)^{\frac{2}{\alpha}}, \quad dq=-\frac{2q_0}{\alpha}(\cos \theta)^{\frac{2-\alpha}{\alpha}} \sin \theta$$
this becomes
\begin{eqnarray*}
\tau&=&\frac{2q_0^{\alpha+1}}{\alpha q_0^{\frac{\alpha}{2}}}\int_0^{\frac{\pi}{2}}
\frac{(\cos \theta)^{\frac{2}{\alpha}} \sin \theta}{\sqrt{1-\cos^2 \theta}}d\theta\\
&=&\frac{2q_0^{\frac{\alpha+2}{2}}}{\alpha}\int_0^{\frac{\pi}{2}}
(\cos \theta)^{\frac{2}{\alpha}} d\theta\\
&=&\frac{q_0^{\frac{\alpha+2}{2}}}{\alpha }B\bigg(\frac{2+\alpha}{2\alpha},\frac{1}{2}\bigg)\\
&=&\frac{q_0^{\frac{\alpha+2}{2}}}{\alpha }
\frac{\Gamma\big(\frac{2+\alpha}{2\alpha}\big)\Gamma\big(\frac{1}{2}\big)}{\Gamma\big(\frac{1+\alpha}{\alpha}\big)}
\end{eqnarray*}
Where $B$ is the Betafunction. Similarly we have
\begin{eqnarray*}
\int_0^{\tau} qdt&=&q_0^\frac{\alpha}{2}\int_0^{q_0}\frac{q^{\frac{\alpha+2}{2}}}{\sqrt{q_0^\alpha-q^\alpha}}dq\\
&=&\frac{2q_0^{\frac{\alpha+4}{2}}}{\alpha}\int_0^{\frac{\pi}{2}}
(\cos \theta)^{\frac{4}{\alpha}} d\theta\\
&=&\frac{q_0^{\frac{\alpha+4}{2}}}{\alpha}B\bigg(\frac{4+\alpha}{2\alpha},\frac{1}{2}\bigg)\\
&=&\frac{q_0^{\frac{\alpha+4}{2}}}{\alpha}
\frac{\Gamma\big(\frac{4+\alpha}{2\alpha}\big)\Gamma\big(\frac{1}{2}\big)}{\Gamma\big(\frac{2+\alpha}{\alpha}\big)}
\end{eqnarray*}
implying that
\begin{eqnarray*}
\overline{q}&=&\frac{1}{\tau}\int_0^{\tau}qdt\\
&=&\frac{\Gamma\big(\frac{4+\alpha}{2\alpha}\big)\Gamma\big(\frac{1+\alpha}{\alpha}\big)}{
\Gamma\big(\frac{2+\alpha}{\alpha}\big)\Gamma\big(\frac{2+\alpha}{2\alpha}\big)}q_0
\end{eqnarray*}
This proves the proposition. \hfill $\square$

\bigskip

\begin{flushleft}

%
Urs Frauenfelder,\\
Mathematisches Institut,\\
Universit\"at Augsburg\\
Universit\"atsstra\ss e 14, \\
86159 Augsburg, Germany\\[1ex]
e-mail: \texttt{urs.frauenfelder@math.uni-augsburg.de}\\[2ex]

%

\end{flushleft}
\end{document}